\documentclass[10pt]{article}
\usepackage{amsmath,amssymb,amscd,verbatim,amsthm}
\usepackage{amsmath}
\usepackage{amssymb}

\begin{document}

\begin{center}
{\large \textbf{The Hyers-Ulam stability for nonlinear Volterra integral
equations via a generalized Diaz-Margolis's fixed point theorem}}\\[0.4in]
\textbf{Wei-Shih Du}\footnote{{\small E-mail address:
wsdu@nknucc.nknu.edu.tw, wsdu@mail.nknu.edu.tw\ (W.-S.\ Du); Tel:
+886-7-7172930 ext 6809; Fax: +886-7-6051061.}}\bigskip

Department of Mathematics, National Kaohsiung Normal University, Kaohsiung
82444, Taiwan
\end{center}

\bigskip

\hrule\vspace{0.1cm}\bigskip

\noindent \textbf{Abstract:} In this work, we prove an existence theorem of
the Hyers-Ulam stability for the nonlinear Volterra integral equations which
improves and generalizes Castro-Ramos theorem by using some weak
conditions.\bigskip

\noindent \textbf{2010 Mathematics Subject Classification:} \bigskip 26D10,
34K20, 45D05, 47H10.

\vspace{0.1cm} \hrule\vspace{0.3cm} \bigskip

\noindent {\large \textbf{1. Introduction and preliminaries}}\bigskip

The stability of functional equations was originally raised by Ulam in 1940
and the problem posed by Ulam was the following: "Under what conditions does
there exist an additive mapping near an approximately additive mapping?" (we
refer the reader to [1] for details). In 1941, Hyers [2] gave the first
answer to Ulam's question in the case of Banach spaces. Since then a number
of generalizations in the study of stability of functional equations have
been investigated by several authors; see [3-5] and references therein. A
generalized Banach contraction principle in a complete generalized metric
space proved by Diaz and Margolis [6] has played an important role in the
study of stability of functional equations.\bigskip

\noindent \textbf{Definition 1.1. [6]}\quad Let $X$ be a nonempty set. A
function $p:X\times X\rightarrow \lbrack 0,\infty ]$ is called a \textit{%
generalized} \textit{metric }on $X$ if the following conditions hold:

\begin{enumerate}
\item[$(GM1)$] $p(x,y)=0$ if and only if $x=y$ ;

\item[$(GM2)$] $p(x,y)=p(y,x)$ for all $x$, $y$ $\in X$;

\item[$(GM3)$] $p(x,z)\leq p(x,y)+p(y,z)$ for all $x,y,z\in X$.
\end{enumerate}

\noindent The pair $(X,p)$ is then called a \textit{generalized} \textit{%
metric space}.\bigskip

We remark that the only one difference of the generalized metric from the
usual metric is that the range of the former is permitted to include the
infinity.\bigskip

\noindent \textbf{Theorem 1.1. (Diaz and Margolis [6])}\quad \textit{Let }$%
(X,p)$\textit{\ be a complete generalized metric space and }$T:X\rightarrow
X $\textit{\ be a selfmap on }$X$\textit{. Assume that there exists a
nonnegative real number }$\lambda <1$ \textit{such that }%
\begin{equation*}
p(Tx,Ty)\leq \lambda p(x,y)\text{ \ \textit{for all} }x,y\in X\text{.}
\end{equation*}%
\textit{Denote }$T^{0}=I$\textit{, the identity mapping. Then, for a given
element }$u\in X$\textit{, exactly one of the following assertions is true:}

\begin{enumerate}
\item[\textit{(a)}] $p(T^{n}u,T^{n+1}u)=\infty $ \ \textit{for all} $n\in
\mathbb{N}
\cup \{0\}$,

\item[\textit{(b)}] \textit{there exists a nonnegative integer }$\ell $
\textit{such that }$p(T^{n}u,T^{n+1}u)<\infty $ \ \textit{for all} $n\geq
\ell $.
\end{enumerate}

\textit{Actually, if the assertion (b) holds, then}

\begin{enumerate}
\item[\textit{(b1)}] \textit{the sequence }$\{T^{n}u\}_{n\in
\mathbb{N}
\cup \{0\}}$\textit{\ is convergent to a fixed point }$\hat{y}$ \textit{of} $%
T$.

\item[\textit{(b2)}] $\hat{y}$ \textit{is the unique fixed point of} $T$
\textit{in the set }$\mathcal{S}$\textit{, where}
\begin{equation*}
\mathcal{S}=\{x\in X:p\left( T^{\ell }u,x\right) <\infty \}\text{;}
\end{equation*}

\item[\textit{(b3)}] $p(x,\hat{y})\leq \frac{1}{1-\lambda }p(x,Tx)$ \textit{%
for all }$x\in \mathcal{S}$.
\end{enumerate}

\bigskip

\noindent \textbf{Definition 1.2. [7, 8]}\quad A function $\varphi :$ $%
[0,\infty )\rightarrow $ $[0,1)$ is said to be an $\mathcal{MT}$-\textit{%
function} (or $\mathcal{R}$-\textit{function}) if $\limsup\limits_{s%
\rightarrow t^{+}}\varphi (s)<1$ for all $t\in \lbrack 0,\infty )$. \bigskip

It is obvious that if $\varphi :$ $[0,\infty )\rightarrow $ $[0,1)$ is a
nondecreasing function or a nonincreasing function, then $\varphi $ is an $%
\mathcal{MT}$-function. So the set of $\mathcal{MT}$-functions is a rich
class. In 2012, Du [8] proved the following characterizations of $\mathcal{MT%
}$-functions.\bigskip

\noindent \textbf{Theorem 1.2. [8]}\quad \textit{Let }$\varphi :$\textit{\ }$%
[0,\infty )\rightarrow $\textit{\ }$[0,1)$\textit{\ be a function. Then the
following statements are equivalent.}

\begin{enumerate}
\item[(a)] $\varphi $\textit{\ is an }$\mathcal{MT}$\textit{-function.}

\item[(b)] \textit{For each }$t\in \lbrack 0,\infty )$\textit{, there exist }%
$r_{t}^{(1)}\in \lbrack 0,1)$\textit{\ and }$\varepsilon _{t}^{(1)}>0$%
\textit{\ such that }$\varphi (s)\leq r_{t}^{(1)}$\textit{\ for all }$s\in
(t,t+\varepsilon _{t}^{(1)})$\textit{.}

\item[(c)] \textit{For each }$t\in \lbrack 0,\infty )$\textit{, there exist }%
$r_{t}^{(2)}\in \lbrack 0,1)$\textit{\ and }$\varepsilon _{t}^{(2)}>0$%
\textit{\ such that }$\varphi (s)\leq r_{t}^{(2)}$\textit{\ for all }$s\in
\lbrack t,t+\varepsilon _{t}^{(2)}]$\textit{.}

\item[(d)] \textit{For each }$t\in \lbrack 0,\infty )$\textit{, there exist }%
$r_{t}^{(3)}\in \lbrack 0,1)$\textit{\ and }$\varepsilon _{t}^{(3)}>0$%
\textit{\ such that }$\varphi (s)\leq r_{t}^{(3)}$\textit{\ for all }$s\in
(t,t+\varepsilon _{t}^{(3)}]$\textit{.}

\item[(e)] \textit{For each }$t\in \lbrack 0,\infty )$\textit{, there exist }%
$r_{t}^{(4)}\in \lbrack 0,1)$\textit{\ and }$\varepsilon _{t}^{(4)}>0$%
\textit{\ such that }$\varphi (s)\leq r_{t}^{(4)}$\textit{\ for all }$s\in
\lbrack t,t+\varepsilon _{t}^{(4)})$\textit{.}

\item[(f)] \textit{For any nonincreasing sequence }$\{x_{n}\}_{n\in
\mathbb{N}
}$\textit{\ in }$[0,\infty )$\textit{, we have }$0\leq \sup\limits_{n\in
\mathbb{N}
}\varphi (x_{n})<1$\textit{.}

\item[(g)] $\varphi $\textit{\ is a function of contractive factor; that is,
for any strictly decreasing sequence }$\{x_{n}\}_{n\in
\mathbb{N}
}$\textit{\ in }$[0,\infty )$\textit{, we have }$0\leq \sup\limits_{n\in
\mathbb{N}
}\varphi (x_{n})<1$\textit{.}
\end{enumerate}

\bigskip

In 2009, Castro and Ramos [9] proved the following existence theorem of the
Hyers-Ulam stability for the nonlinear Volterra integral equations.\bigskip

\noindent \textbf{Theorem 1.3. [9, Theorem 5.1]}\quad \textit{Let }$a$%
\textit{\ and }$b $\textit{\ be given real numbers with }$a<b$\textit{\ and
let} $K:=b-a$. \textit{Assume that there exists a positive constant }$L$
\textit{such that} \textit{\ }%
\begin{equation*}
0<KL<1\text{\textit{.}}
\end{equation*}%
\textit{Assume that }$f:[a,b]\times \lbrack a,b]\times
\mathbb{C}
\rightarrow
\mathbb{C}
$\textit{\ is a continuous function satisfying }%
\begin{equation*}
\left\vert f(x,\tau ,y)-f(x,\tau ,z)\right\vert \leq L\left\vert
y-z\right\vert \text{\ \ \textit{for} \textit{any }}x,\tau \in \lbrack a,b]%
\text{ \textit{and} }y,z\in
\mathbb{C}
\text{\textit{.}}
\end{equation*}%
\textit{If there exists a continuous function }$y:[a,b]\rightarrow
\mathbb{C}
$ \textit{satisfying }%
\begin{equation*}
\left\vert y(x)-\int_{a}^{x}f(x,\tau ,y(\tau ))d\tau \right\vert \leq \theta
\text{\ }
\end{equation*}%
\textit{for each }$x\in \lbrack a,b]$\textit{\ and some constant }$\theta
\geq 0$\textit{, then there exists a unique continuous function }$%
y_{0}:[a,b]\rightarrow
\mathbb{C}
$\textit{\ such that }%
\begin{equation*}
y_{0}(x)=\int_{a}^{x}f(x,\tau ,y_{0}(\tau ))d\tau
\end{equation*}%
\textit{and }%
\begin{equation*}
\left\vert y(x)-y_{0}(x)\right\vert \leq \frac{\theta }{1-KL}
\end{equation*}%
\textit{for all }$x\in \lbrack a,b]$\textit{.}\bigskip \bigskip

In this work, we give a generalization of Castro-Ramos theorem by using some
weak conditions.\bigskip \bigskip \bigskip

\noindent {\large \textbf{2. Main results}}\bigskip

Very recently, Du [10] established the following generalization of
Diaz-Margolis's fixed point theorem.\bigskip

\noindent \textbf{Theorem 2.1. [10]}\quad \textit{Let }$(X,p)$\textit{\ be a
complete generalized metric space and }$T:X\rightarrow X$\textit{\ be a
selfmap on }$X$\textit{. Assume that there exists an} $\mathcal{MT}$-\textit{%
function} $\alpha :[0,\infty )\rightarrow \lbrack 0,1)$ \textit{such that }%
\begin{equation*}
p(Tx,Ty)\leq \alpha (p(x,y))p(x,y)\text{ \ \textit{for all} }x,y\in X\text{
\textit{with} }p(x,y)<\infty \text{.}
\end{equation*}%
\textit{Denote }$T^{0}=I$\textit{, the identity mapping. Then, for a given
element }$u\in X$\textit{, exactly one of the following assertions is true:}

\begin{enumerate}
\item[\textit{(a)}] $p(T^{n}u,T^{n+1}u)=\infty $ \ \textit{for all} $n\in
\mathbb{N}
\cup \{0\}$,

\item[\textit{(b)}] \textit{there exists a nonnegative integer }$\ell $
\textit{such that }$p(T^{n}u,T^{n+1}u)<\infty $ \ \textit{for all} $n\geq
\ell $.
\end{enumerate}

\textit{Actually, if the assertion (b) holds, then}

\begin{enumerate}
\item[\textit{(b1)}] \textit{the sequence }$\{T^{n}u\}_{n\in
\mathbb{N}
\cup \{0\}}$\textit{\ is convergent to a fixed point }$v$ \textit{of} $T$.

\item[\textit{(b2)}] $v$ \textit{is the unique fixed point of} $T$ \textit{%
in the set }$\mathcal{L}$\textit{, where}
\begin{equation*}
\mathcal{L}=\{x\in X:p(T^{\ell }u,x)<\infty \}\text{;}
\end{equation*}

\item[\textit{(b3)}] $p(x,v)\leq \frac{1}{1-\alpha (p(x,v))}p(x,Tx)$ \textit{%
for all }$x\in \mathcal{L}$.
\end{enumerate}

\bigskip

In this paper, we prove an existence theorem of the Hyers-Ulam stability for
the nonlinear Volterra integral equations which improves and generalizes [9,
Theorem 5.1].\bigskip

\noindent \textbf{Theorem 2.2.}\quad \textit{Let }$a$\textit{\ and }$b$%
\textit{\ be given real numbers with }$a<b$\textit{\ and let} $K:=b-a$.
\textit{Assume that }$\varphi :%
\mathbb{R}
\rightarrow
\mathbb{R}
$\textit{\ is a function which\ is nondecreasing on }$[0,\infty )$\textit{\
satisfying} \
\begin{equation*}
\varphi \left( \lbrack 0,\infty )\right) \subseteq \left[ 0,\delta K^{-1}%
\right] \text{\textit{,}}
\end{equation*}%
\textit{for some constant }$0<\delta <1$, \textit{and }$V:[a,b]\times
\lbrack a,b]\times
\mathbb{C}
\rightarrow
\mathbb{C}
$\textit{\ is a continuous function satisfying }%
\begin{equation*}
\left\vert V(x,\tau ,y)-V(x,\tau ,z)\right\vert \leq \varphi \left(
\left\vert y-z\right\vert \right) \left\vert y-z\right\vert \text{\ \
\textit{for} \textit{any }}x,\tau \in \lbrack a,b]\text{ \textit{and} }%
y,z\in
\mathbb{C}
\text{\textit{.}}
\end{equation*}%
\textit{If there exists a continuous function }$y:[a,b]\rightarrow
\mathbb{C}
$ \textit{satisfying }%
\begin{equation}
\left\vert y(x)-\int_{a}^{x}V(x,\tau ,y(\tau ))d\tau \right\vert \leq \theta
\text{\ }  \tag{2.1}
\end{equation}%
\textit{for each }$x\in \lbrack a,b]$\textit{\ and some constant }$\theta
\geq 0$\textit{, then there exists a unique continuous function }$%
y_{0}:[a,b]\rightarrow
\mathbb{C}
$\textit{\ such that }%
\begin{equation*}
y_{0}(x)=\int_{a}^{x}V(x,\tau ,y_{0}(\tau ))d\tau
\end{equation*}%
\textit{and }%
\begin{equation*}
\left\vert y(x)-y_{0}(x)\right\vert \leq \frac{\theta }{1-\delta }
\end{equation*}%
\textit{for all }$x\in \lbrack a,b]$\textit{.}\medskip

\noindent \textbf{Proof.}\quad Let $X$ denote the set of all continuous
functions from $[a,b]$ to $%
\mathbb{C}
$. Define a function $p:X\times X\rightarrow \lbrack 0,\infty ]$ by
\begin{equation*}
p(f,g)=\inf \{M\geq 0:\left\vert f(x)-g(x)\right\vert \leq M\text{ for all }%
x\in \lbrack a,b]\},
\end{equation*}%
where we adopt the usual convention that $\inf \emptyset =\infty $. Then $%
(X,p)$\ is a complete generalized metric spacespace. Let $T:X\rightarrow X$
be defined by
\begin{equation*}
(Tf)(x)=\int_{a}^{x}V(x,\tau ,f(\tau ))d\tau
\end{equation*}%
for all $f\in X$ and $x\in \lbrack a,b]$. It is esay to show $Tf\in X$ for
all $f\in X\mathit{\ }$and this ensures that $T$ is well defined. Define an $%
\mathcal{MT}$-function $\alpha :[0,\infty )\rightarrow \lbrack 0,1)$ by
\begin{equation*}
\alpha (t)=K\varphi (t)\text{ \ for all }t\in \lbrack 0,\infty )\text{.}
\end{equation*}%
It is not hard to verify that for all $f,g\in X$ with $p(f,g)<\infty $,
\textit{\ }%
\begin{equation*}
p(Tf,Tg)\leq \alpha (p(f,g))p(f,g)\text{.}
\end{equation*}%
Take $h\in X$. Since $p(Tf,f)<\infty $ for all $f\in X$, we have $%
p(Th,h)<\infty $. For any $f\in X$, since $f$ and $h$ are continuous on $%
[a,b]$, there exists a constant $c\geq 0$ such that
\begin{equation*}
\left\vert h(x)-f(x)\right\vert \leq c\text{ \ for any }x\in \lbrack a,b]
\end{equation*}%
which implies $p(h,f)\leq c<\infty $. Therefore we prove
\begin{equation*}
X=\{f\in X:p(h,f)<\infty \}\text{.}
\end{equation*}%
Applying Theorem 2.1 (b), there exists a unique\textit{\ }$y_{0}\in X$ such
that $T^{n}h\overset{p}{\longrightarrow }y_{0}$ \ \ as $n\rightarrow \infty $%
, $Ty_{0}=y_{0}$ and
\begin{equation}
p(f,y_{0})\leq \frac{1}{1-\alpha (p(f,y_{0}))}p(f,Tf)\text{ \ \ for all }%
f\in X.  \tag{2.2}
\end{equation}%
\bigskip So
\begin{equation*}
y_{0}(x)=\int_{a}^{x}V(x,\tau ,y_{0}(\tau ))d\tau \text{ \ \ for all\textit{%
\ }}x\in \lbrack a,b]\text{.}
\end{equation*}%
By (2.1), we have \textit{\ }%
\begin{equation}
p(y,Ty)\leq \theta .  \tag{2.3}
\end{equation}%
Since $\alpha (p(y,y_{0}))\leq \delta $, by (2.2) and (2.3), we get
\begin{equation*}
p(y,y_{0})\leq \frac{1}{1-\alpha (p(y,y_{0}))}p(y,Ty)\leq \frac{\theta }{%
1-\delta }\text{,}
\end{equation*}%
which deduce
\begin{equation*}
\left\vert y(x)-y_{0}(x)\right\vert \leq \frac{\theta }{1-\delta }\text{ \
for all\textit{\ }}x\in \lbrack a,b]\text{.}
\end{equation*}%
The proof is completed. \hspace{\fill}\textit{\ } $\Box $\bigskip

\noindent \textbf{Remark 2.3.}\quad \lbrack 9, Theorem 5.1] is a special
case of Theorem 2.2. Indeed, let $g:%
\mathbb{R}
\rightarrow
\mathbb{R}
$ be any function. Define $\varphi :%
\mathbb{R}
\rightarrow
\mathbb{R}
$\ by
\begin{equation*}
\varphi (t)=\left\{
\begin{array}{cc}
L, & \text{for }t\geq 0, \\
g(t)\text{,} & \text{otherwise.}%
\end{array}%
\right. .
\end{equation*}%
Put $\delta :=KL$. Then all the conditions of Theorem 2.2 are satisfied and
the conclusion of [9, Theorem 5.1] follows from Theorem 2.2.\bigskip
\bigskip \bigskip

\noindent {\large \textbf{References}}

\begin{enumerate}
\item[{[1]}] S.M. Ulam, A Collection of Mathematical Problems, Interscience
Tracts in Pure and Applied Mathematics 8, Interscience Publishers, New York,
1960.

\item[{[2]}] D.H. Hyers, On the stability of the linear functional equation,
Proc. Natl. Acad. Sci. USA 27 (1941) 222-224.

\item[{[3]}] V. Lakshmikantham and M.R.M. Rao, Theory of Integro-differential
Equations, Stability and Control: Theory, Methods and Applications 1, Gordon
and Breach Publ., Philadelphia, 1995.

\item[{[4]}] S.-M. Jung, Hyers-Ulam-Rassias Stability of Functional Equations
in Nonlinear Analysis, in: Springer Optimization and Its Applications, vol.
48, Springer, New York, 2011.

\item[{[5]}] I.A. Rus, Ulam stability of the operatorial equations, in:
Functional Equations in Mathematical Analysis, in: Springer Optim. Appl.,
vol. 52, Springer, New York, 2012, pp. 287-305.

\item[{[6]}] J.B. Diaz and B. Margolis, A fixed point theorem of the
alternative, for contractions on a generalized complete metric space, Bull.
Amer. Math. Soc. 74 (1968) 305-309.

\item[{[7]}] W.-S. Du, Some new results and generalizations in metric fixed
point theory, Nonlinear Anal. 73 (2010) 1439-1446.

\item[{[8]}] W.-S. Du, On coincidence point and fixed point theorems for
nonlinear multivalued maps, Topology and its Applications 159 (2012) 49-56.

\item[{[9]}] L.P. Castro and A. Ramos, Hyers-Ulam-Rassias stability for a
class of nonlinear Volterra integral equations, Banach J. Math. Anal.(3)
(2009), no. 1, 36-43.

\item[{[10]}] W.-S. Du, The generalization of Diaz-Margolis's fixed point
theorem and its application to the stability for generalized Volterra
integral equations, submitted.
\end{enumerate}

\end{document}